%% file: author.tex
\documentclass[graybox]{svmult}

\usepackage{mathptmx}       
\usepackage{helvet}         
\usepackage{courier}        
\usepackage{type1cm}        
\usepackage{makeidx}         
\usepackage{graphicx}        
\usepackage{multicol}        
\usepackage[bottom]{footmisc}

\usepackage{amsmath}
\usepackage{amssymb}
\usepackage{amsfonts}
\usepackage{dsfont}
\usepackage{subfigure}

\include{macros}

\makeindex             

\begin{document}

\title*{A Production Model with History Based Random Machine Failures}
\author{Stephan Knapp and Simone G\"ottlich}
\institute{S.\ Knapp \at University of Mannheim, 68131 Mannheim, \email{stknapp@mail.uni-mannheim.de}
\and S.\ G\"ottlich \at University of Mannheim, 68131 Mannheim \email{goettlich@uni-mannheim.de}}
%
%
\maketitle

\abstract*{In this paper, we introduce a time-continuous production model that enables random machine failures, where the failure probability depends historically on the production itself. This bidirectional relationship between historical failure probabilities and production is mathematically modeled by the theory of piecewise deterministic Markov processes (PDMPs). On this way, the system is rewritten into a Markovian system  such that classical results can be applied. In addition, we present a suitable solution, taken from machine reliability theory, to connect past production and the failure rate. Finally, we investigate the behavior of the presented model numerically in examples by considering sample means of relevant quantities and relative frequencies of number of repairs.}

\abstract{In this paper, we introduce a time-continuous production model that enables random machine failures, where the failure probability depends historically on the production itself. This bidirectional relationship between historical failure probabilities and production is mathematically modeled by the theory of piecewise deterministic Markov processes (PDMPs). On this way, the system is rewritten into a Markovian system  such that classical results can be applied. In addition, we present a suitable solution, taken from machine reliability theory, to connect past production and the failure rate. Finally, we investigate the behavior of the presented model numerically in examples by considering sample means of relevant quantities and relative frequencies of number of repairs.}

\section{Modeling Equations} \label{sec:The Model}
We briefly recall the production network model from \cite{ApiceGoettlichHertyBenedetto,GoettlichHertyKlar2005} first, and according to \cite{GoettlichKnapp2018}, we present the stochastic extension to a load-dependent production model with machine failures. To keep the notation well-arranged, we consider a production network consisting of a single queue processor unit. 
We assume a processor, which is represented by an interval $(a,b) \subset \R$, i.e., with length $L = b-a$, where $\rho(x,t)$ describes the density of production goods at $x \in (a,b)$ and time $t \geq 0$. The dynamics of the density, and consequently of the production, is given by the following nonlinear hyperbolic partial differential equation
\begin{equation}
\partial_t \rho(x,t) + \partial_x \min\{v \rho(x,t),\mu(t)\}=0, \label{eq:rho}
\end{equation}
where $\mu(t) \geq 0$ is the time-dependent bounded production capacity and $v >0$ the constant production velocity.
In front of the processor a storage, also called queue, is assumed and for an externally given time-dependent inflow $G_{\text{in}}(t)$ into the production, the queue length $q$ follows the ordinary differential equation
\begin{equation}
\partial_t q(t) = G_{\text{in}}(t)-g_{\text{out}}(t), \label{eq:q}
\end{equation}
with
\begin{equation*}
g_{\text{out}}(t) = 
\begin{cases}
\min\{G_{\text{in}}(t),\mu(t)\}, &\text{ if } q(t) = 0,\\
\mu(t), &\text{ if } q(t)>0.
\end{cases}
\end{equation*}
The processor is coupled to the queue by a boundary condition in the form of $\rho(a,t) = \frac{g_{\text{out}}(t)}{v}$ and initial conditions $\rho(x,0) = \rho_0(x) \in L^1((a,b))$, $q(0) = q_0 \in \R_{\geq 0}$ are prescribed.
This model is well-defined, see, e.g.\ \cite{ApiceGoettlichHertyBenedetto}, if the capacity is independent of time.
The theory of piecewise deterministic Markov processes; see, e.g. \cite{Davis1984,Jacobsen2006}, has been used to define an appropriate production model with stochastic machine failures in \cite{GoettlichKnapp2018}, where the probabilities of machine failures depend on the actual workload of the processor. 
Since this construction only allows for a dependence on the current workload, we can not use the amount of goods produced since the last machine failure as a measure for the next failure. Our new idea lies in adding a variable $w$ governing the workload since the last repair. To do so, we use the time-dependent variable $r(t) \in \{0,1\}$, and set the capacity as $\mu(t) = r(t) c$ for a maximal capacity $c>0$. This means that $r(t) = 0 \Rightarrow \mu(t) = 0$ is a down and $r(t) = 1 \Rightarrow \mu(t) = c$ a working processor at time $t$ and
we define
\begin{equation*}
\WIP(t_0,t_1) = \int_{t_0}^{t_1} \int_a^b \rho(x,t)dxdt
\end{equation*}
as the cumulative work-in-progress of the processor between time $t_0$ and $t_1$. The variable $w$ should therefore satisfy
\begin{align}
\partial_t w(t) = r(t)\int_a^b \rho(x,t)dx, \quad w(t_0) = w_0 = \int_a^b \rho(x,t_0)dx.\label{eq:w}
\end{align}
Altogether, we define the state space
\begin{equation*}
E = \R_{\geq 0}\times \{0,1\} \times \R_{\geq 0} \times L^1((a,b)),
\end{equation*}
which is a measurable space together with the $\sigma$-algebra $\cE$ generated by the open sets induced by the metric 
$$d((w,r,q,\rho),(\tilde{w},\tilde{r},\tilde{q},\tilde{\rho})) = |w-\tilde{w}|+|r-\tilde{r}|+|q-\tilde{q}|+\|\rho-\tilde{\rho}\|_{L^1((a,b))}.$$ 

Since we construct a piecewise deterministic Markov process, we define the deterministic dynamics between jump times as
\begin{align*}
\Phi_{st} \colon E &\to E,\\
(w_0,r_0,q_0,\rho_0) & \mapsto (w(t),r(t),q(t),\rho(t)),
\end{align*}
i.e., $\Phi_{st}$ is the solution to equations \eqref{eq:rho}, \eqref{eq:q}, \eqref{eq:w}, and $ r(t) = r_0$ with initial conditions $(w_0,r_0,q_0,\rho_0) \in E$. 
To characterize the stochastic part, we introduce
\begin{align*}
\psi(t,y) = \gl_{r,r}(t,w),\quad \eta(t,y,B) = \frac{\gl_{r,(1-r)}(t,w)}{\psi(t,y)} \epsilon_{(rw,(1-r),q,\rho)}(B)
\end{align*}
for every $y = (w,r,q,\rho) \in E$ and $B \in \cE$, where $\gl_{i,j}(t,w)$ describes the transition rate from capacity $i$ to $j$ at time $t$ and actual workload $w$, $i,j \in \{0,1\}$ and $\epsilon_x$ is the Dirac measure with unit mass in $x$.
The function $\psi$ is the total intensity determining whether a jump occurs, or not, and the function $\eta$ describes the probability distribution of the systems jump given the system changes at time $t$. For example, given the state $y = (w,1,q,\rho)$ at the time of a jump, the system jumps to $(w,0,q,\rho)$ and, vice versa, given the state $y = (w,0,q,\rho)$ the system jumps to $(0,1,q,\rho)$, i.e., the workload has been ``reset''. The open question is whether this model can be represented by a piecewise deterministic Markov process.

Following \cite{GoettlichKnapp2018}, it is straightforward to show
\begin{theorem}
Let $\gl_{i,j} \colon [0,T] \times \R_{\geq 0} \to \R_{\geq 0}$ be uniformly bounded, continuous and satisfy $\gl_{i,i} = \gl_{i,i-1}$ for $i \in \{0,1\}$. Then for all initial data $x_0 \in E$ there exists a Markov process $$
X = ((w(t),r(t),q(t),\rho(r)),t \in [0,T]) \subset E$$ on some probability space \OAP, satisfying
\begin{enumerate}
\item $X(0) = x_0\quad P$-almost surely,
\item for every $t \in (0,T), (w,r,q,\rho) \in E$ and $j \in \{0,1\}$, it holds that
\begin{align}
P(r(t+\Delta t) = j | X(t) = (w,r,q,\rho)) =&\; \big(1-\Delta t \gl_{r,r}(t,w)\big)\Ind_{r}(j) \notag \\ &\;+\Delta t \gl_{r,(1-r)}(t,w)\Ind_{1-r}(j) + \lano(\Delta t),\notag
\end{align}
\item there exists a $P$-null set $\cN \in \cA$ such that for every $\go \in \gO \setminus \cN$, there exist times $T_0 = 0\leq T_1 \leq \cdots \leq  T_{M}=T$ such that for every $k = 0,\dots,M-1$, 
$X(t) = \Phi_{T_k,t}(X(T_k))$
for $t \in [T_k,T_{k+1})$ with capacity $\mu(r(T_k,\go))$, i.e., $X$ behaves deterministic between jump times.
\end{enumerate}
\end{theorem}
The main and new ingredient is the mapping $t\mapsto w(t)$, which is a continuous mapping since $t \mapsto \rho(t)$ is continuous.

\section{Computational Results} \label{sec:CompRes}
Due to the fact that solutions to \eqref{eq:rho} move with non-negative velocities only, we can use the first order left-sided upwind scheme for a numerical approximation of the density $\rho$. Furthermore, we use the explicit Euler scheme to approximate the queue length $q$ given by \eqref{eq:q} and $w$ given by \eqref{eq:w}, where we use a rectangular rule for the integration. This yields an approximation of the deterministic dynamics between the jump times. The simulation of the jump times is done with the thinning algorithm presented in \cite{GoettlichKnapp2018}. Its basic idea is to use the uniform bound on the rate functions and generate exponentially distributed times with high intensity, representing the times between jumps, and thin these times during the numerical simulation of the whole system with an appropriate acceptance rejection procedure.

The choice of the rate functions $\lambda_{i,j}(t,w)$ is a crucial point in numerical examples. Here, we make use of the choice in \cite{Rivera2018} and set for $\theta_1,\theta_2 >0$ the rate function as
\begin{align*}
\lambda_{1,0}(t,w) = \lambda_{1,0}^{\text{min}}+(\lambda_{1,0}^{\text{max}}-\lambda_{1,0}^{\text{min}})(1-e^{-(\theta_1 w)^{\theta_2}}),
\end{align*}
which is a scaled version of the cumulative distribution function of a Weibull distribution, i.e., $F(t) = 1-e^{-(\theta_1 t)^{\theta_2}}$.
The classical interpretation of $t$ in the latter expression is the lifetime of a machine and $F(t)$ is the probability that a failure happens after time $t$, see, e.g.\ \cite{Jiang2011}. In our case we use the variable $w$, which measures the amount of goods produced since the last repair happened. Therefore, if $w = 0$, then $\lambda_{1,0}(t,0) = \lambda_{1,0}^{\text{min}}$, which corresponds to the minimal failure rate and $\lim_{w \to \infty }\lambda_{1,0}(t,w) = \lambda_{1,0}^{\text{max}}$. The function $\lambda_{1,0}(t,w) $ is monotonically increasing in $w$ and incorporates the idea of an increasing failure rate depending on past workloads.
On the other hand, we assume $\lambda_{0,1}(t,w) = \lambda_{0,1}$ because repair times do not dependent on the amount of goods produced.

In the following, we examine the presented model using numerical examples. Here, we assume a production velocity of $v = 1$, the interval $a = 0$, $b = 1$, and the capacity is given as $\mu(t) = 2 r(t)$. We use a spatial discretization with step-size $\Delta x = 10^{-1}$ and a temporal step-size that satisfies the Courant-Friedrichs-Lewy condition, which reads as $\Delta t \leq \Delta x$ for the chosen parameters. 
The simulation results are based on samples of the stochastic process $X$ and we use the classical Monte-Carlo estimator to evaluate moments or probabilities of the samples. We used a sample size of $10^5$ for all following results. 

\begin{figure}[htb!]
\subfigure[Expected $w$]{
\includegraphics[width=.5\textwidth]{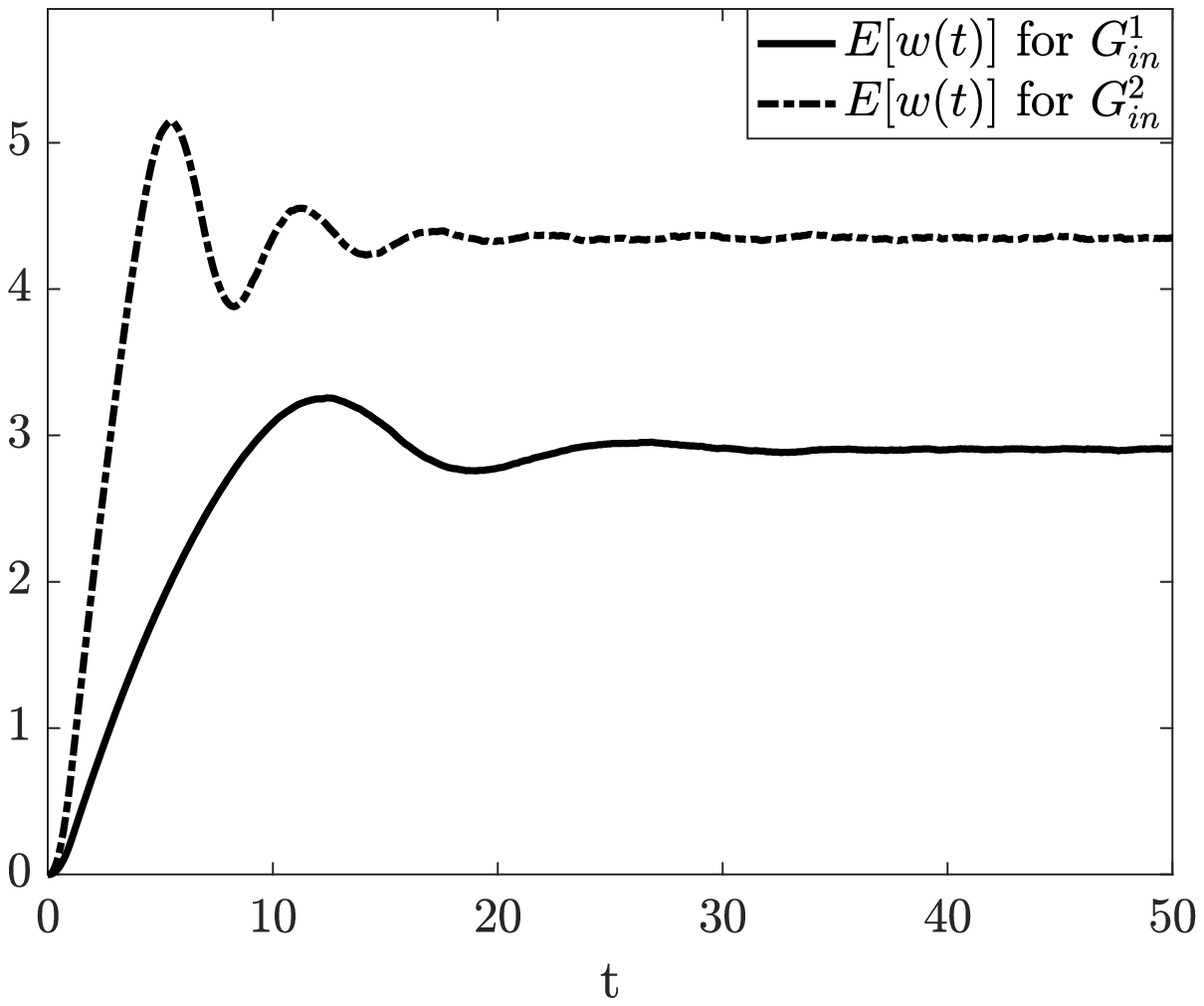}
}
\subfigure[Expected capacity $\mu(t)$]{
\includegraphics[width=.5\textwidth]{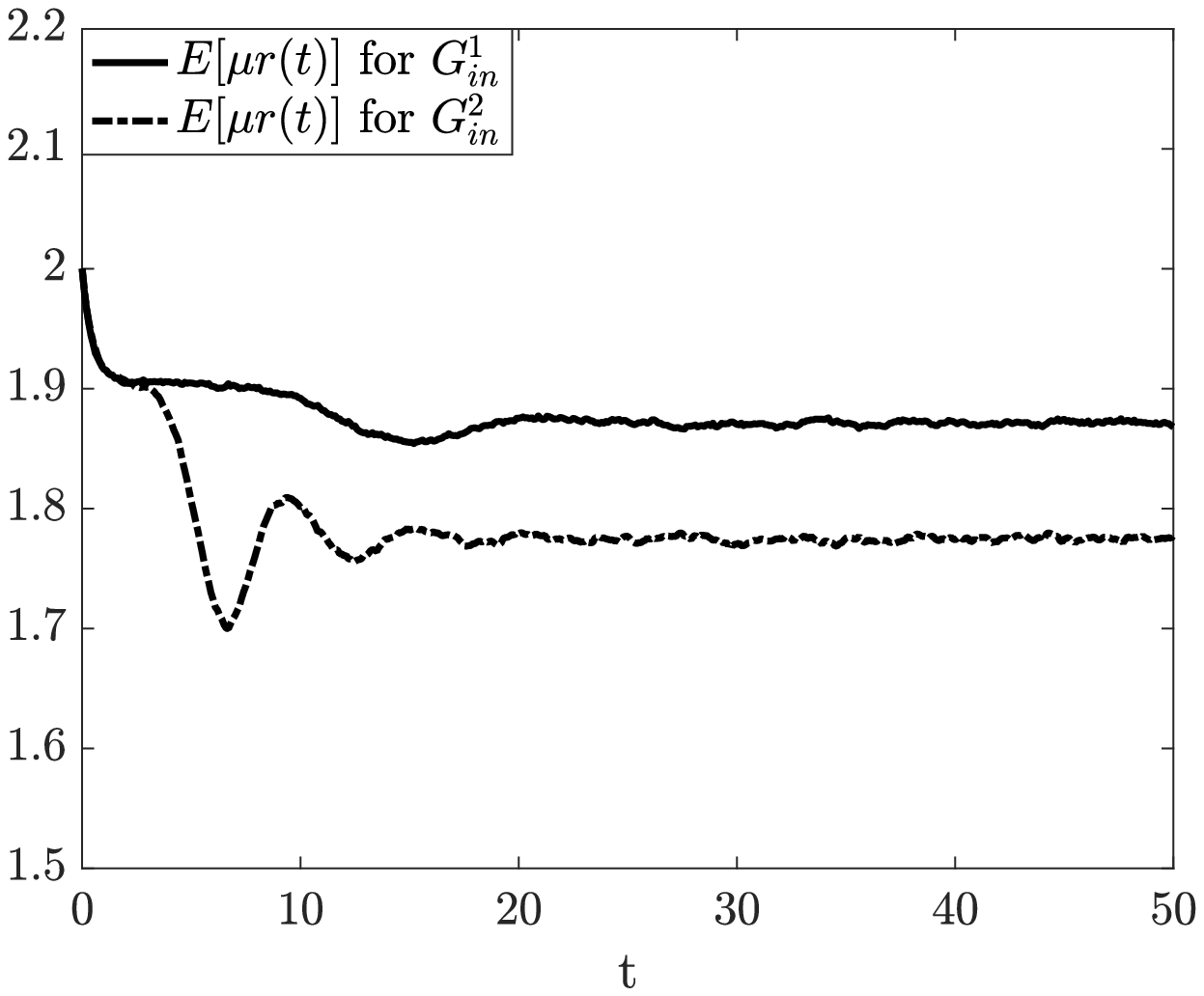}
}
\subfigure[Expected queue-length $q(t)$]{
\includegraphics[width=.5\textwidth]{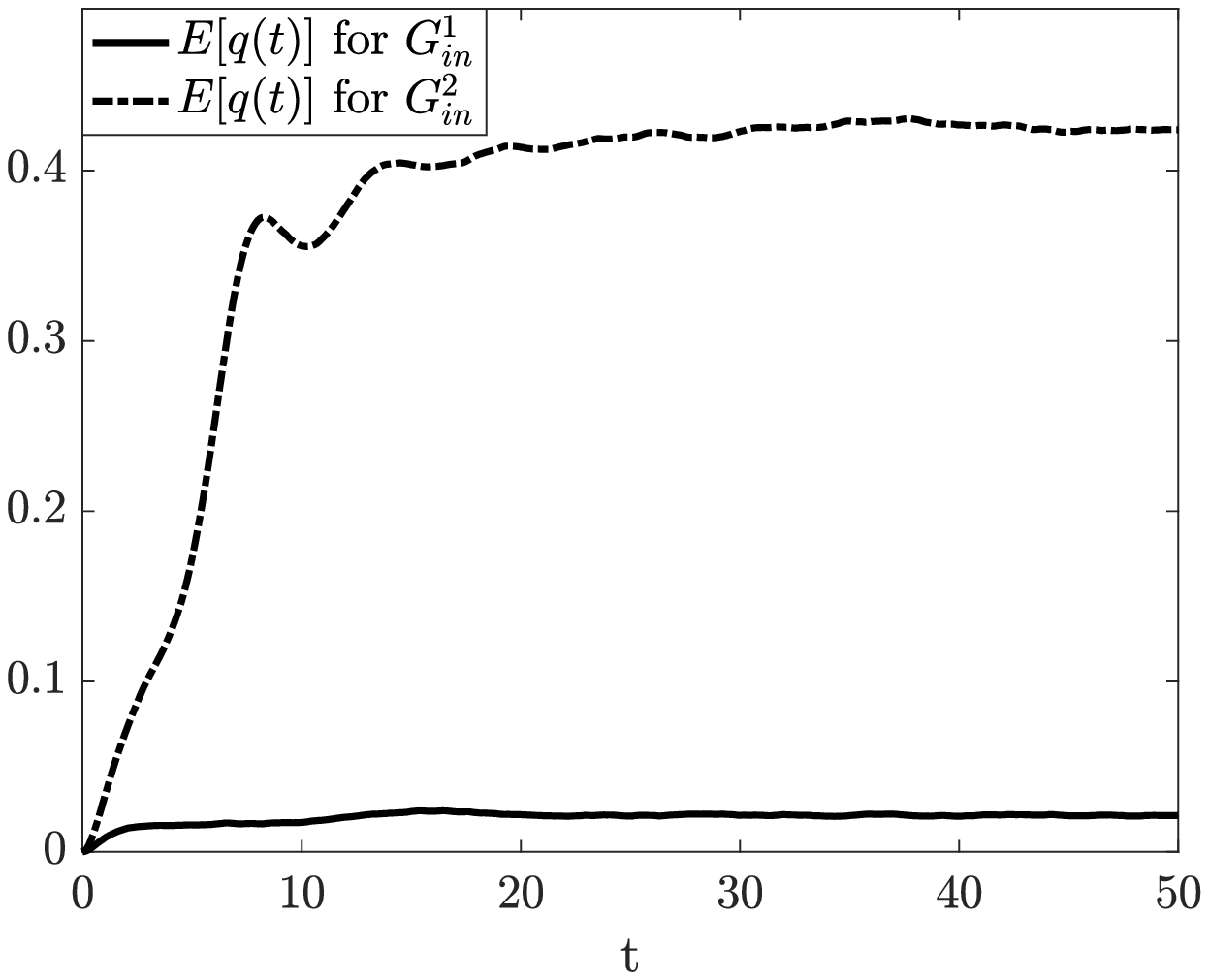}
}
\subfigure[Expected density at $x = 1$]{
\includegraphics[width=.5\textwidth]{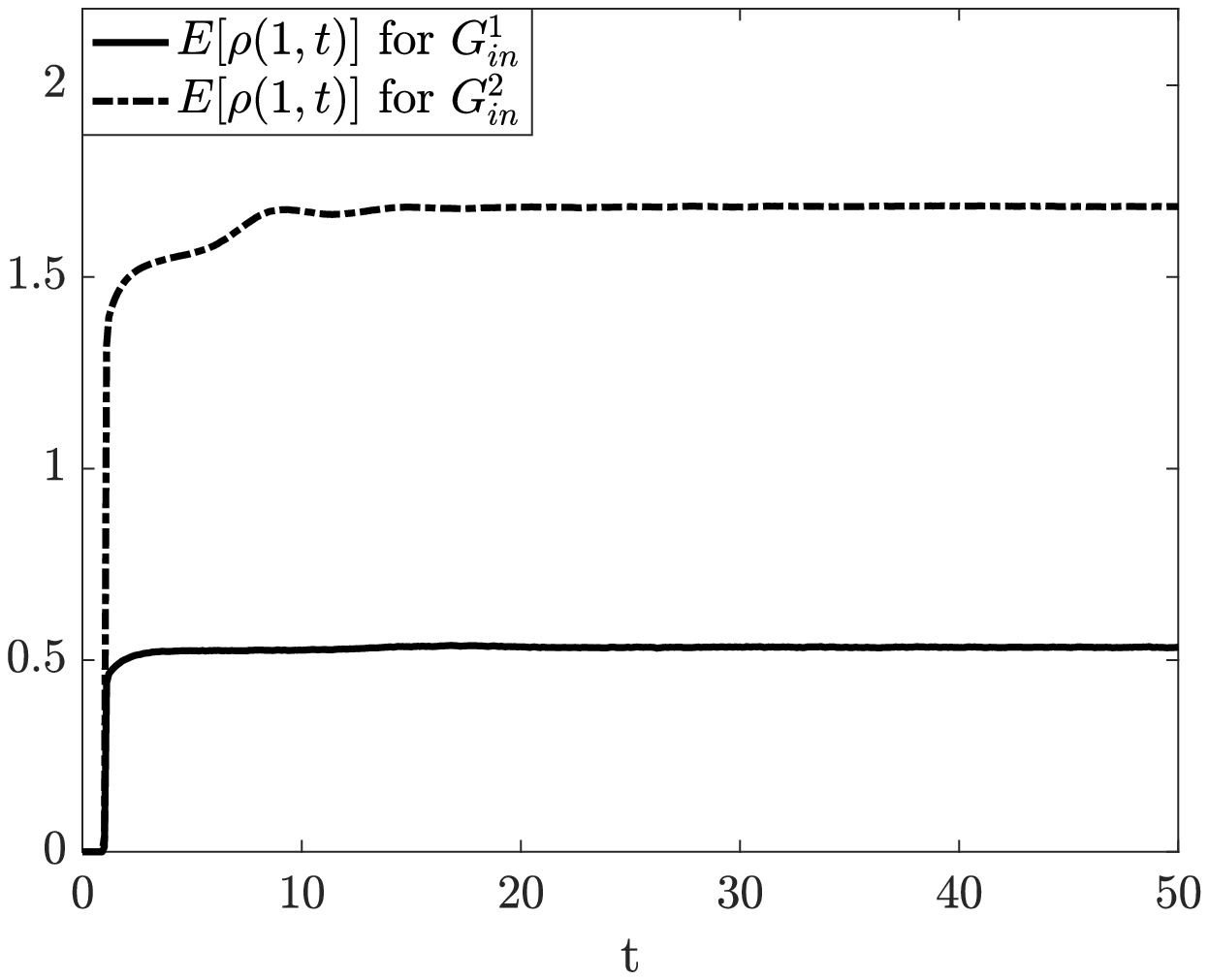}
}
\caption{First order moments of $w$, the capacity, queue-length and density}
\label{fig:ExpectedQuantities1}
\end{figure}

We analyze the expected queue length, capacity and the distribution of the number of repairs within a time horizon $[0,50]$ for two different constant inflow profiles. We denote by $G^{1}_{\text{in}}(t) \equiv 0.5$ and by $G^2_{\text{in}}(t) \equiv 1.5$ as inflow profiles and use the parameters
\begin{align*}
\lambda_{0,1}(t,w) = \frac{1}{0.5},\quad \lambda_{1,0}^{\text{min}} = \frac{1}{10}, \quad \lambda_{1,0}^{\text{max}} = \frac{1}{0.5}, \quad \theta_1 = \frac{1}{10}, \quad \theta_2 = 5.
\end{align*}
In figure \ref{fig:ExpectedQuantities1}, first order moment estimations are shown. In detail, figure \ref{fig:ExpectedQuantities1} (a) shows the expected value of the variable $w$, \ref{fig:ExpectedQuantities1} (b) the expected capacity, \ref{fig:ExpectedQuantities1} (c) the expected queue length and \ref{fig:ExpectedQuantities1} (d) the expected density at the end of the processor. The dynamics is quite interesting: the expected capacity decreases approximately until time $t = 6$ for the second inflow, then increases and decreases again. Indeed, the mean time to failure is given by $\Gamma(1+\frac{1}{\theta_2})\theta_1^{-1}$, see e.g.\ \cite{Jiang2011}. If $w$ corresponds to the lifetime in our model, we see that an intact system with constant inflow $G_{\text{in}}$ is more likely to fail around time $\Gamma(1+\frac{1}{\theta_2})(\theta_1 G_{\text{in}})^{-1}$. In our case, this leads to time $18.4$ for the first and time $6.1$ for the second inflow profile, which is close to the times at which the shape of the expected capacity changes. We observe these characteristic times also in the other graphs in figure \ref{fig:ExpectedQuantities1}. In contrast to the models presented in \cite{ GoettlichKnapp2017,GoettlichKnapp2018, GoettlichMartinSickenberger}, where quantities monotonically converge, we obtain an oscillatory behavior of the quantities for constant inputs. The oscillatory effects are natural and caused by the history we incorporate in $w$. This means, the first machine failures are likely around time $18.4 (6.1)$, the second around $36.8 (12.2)$ and so on. At the same time the failures, which occur between these likely times, smooth this effect out as time evolves and the quantities converge. 

\begin{figure}[htb!]
\subfigure[Inflow $G^1_{\text{in}}$]{
\includegraphics[width=.5\textwidth]{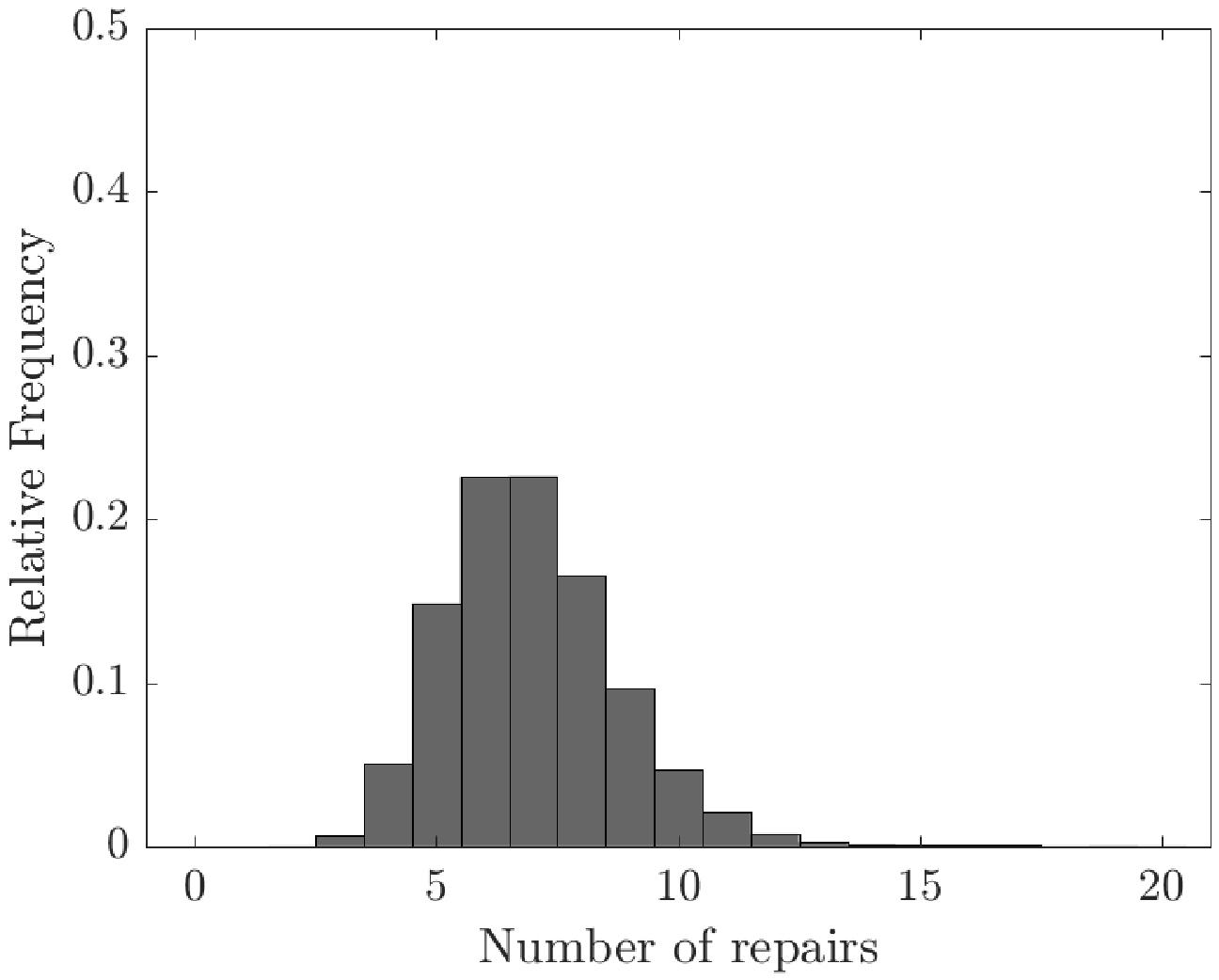}
}
\subfigure[Inflow $G^2_{\text{in}}$]{
\includegraphics[width=.5\textwidth]{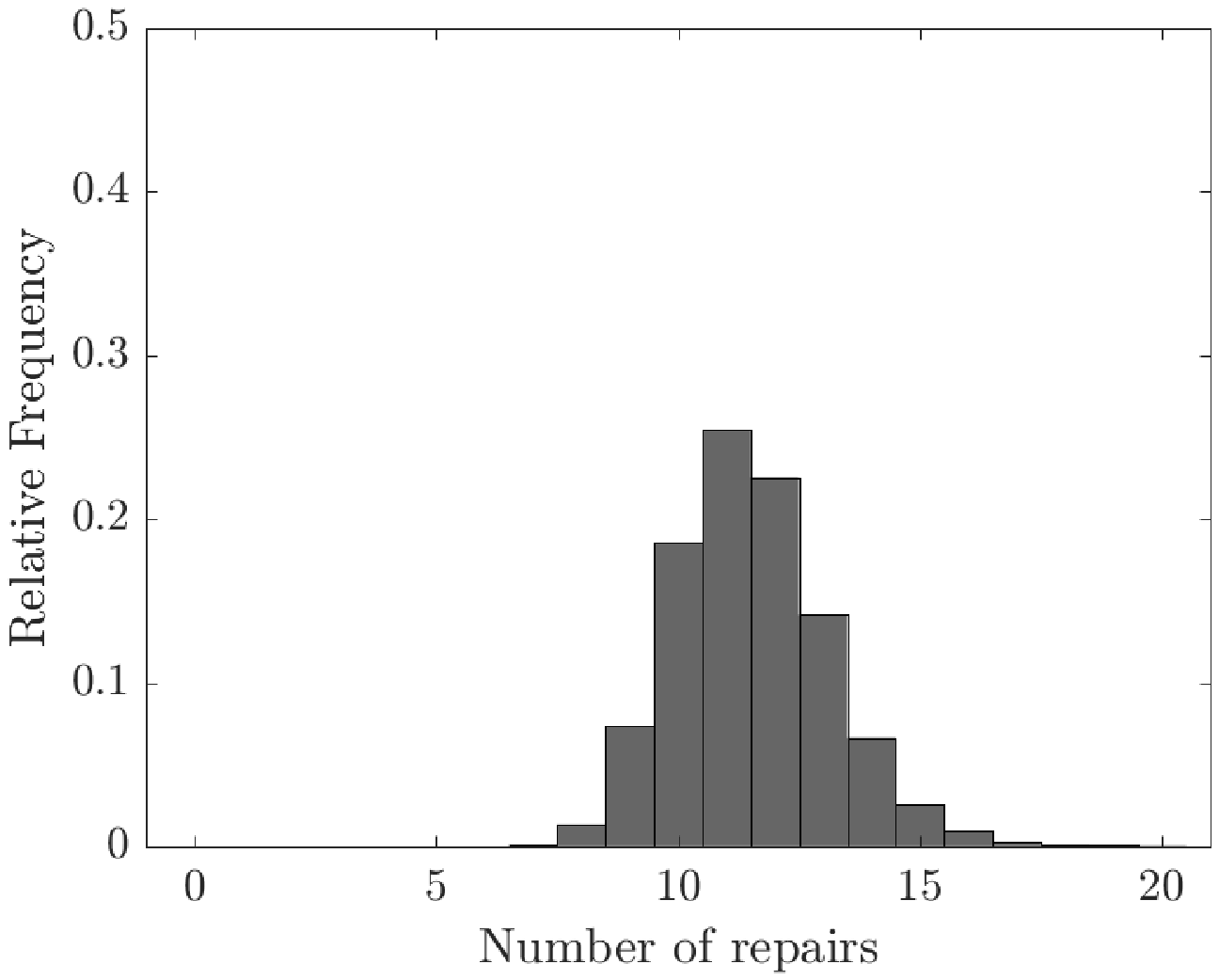}
}
\caption{Distribution of the number of repairs within $[0,50]$}
\label{fig:Histogram1}
\end{figure}

Figure \ref{fig:Histogram1} shows the distribution of the number of repairs within the time horizon $[0,50]$ and emphasizes the impact of the chosen inflow on the reliability of the processor. In figure \ref{fig:Histogram1} (a) the case of $G^1_{\text{in}}$ is shown, where mostly 5 to 9 repairs have been done. The situation for inflow profile $G^2_{\text{in}}$ is different, where 9 to 14 repairs during the time horizon are more likely.

To conclude, we deduced a production model with random machine failures including failure probabilities depending on the workload of the machine since the last repair occured. The extension of the model to complex production networks is straightforward, see, e.g. \cite{GoettlichKnapp2018}. Simulation results showed a big impact of the history on expected workload, capacity, queue length and density. These effects are not negligible for production planning and control and must be taken into account.

\begin{acknowledgement}
This work has been financially supported by the BMBF project ENets (05M18VMA).
\end{acknowledgement}

\input{referenc}
\end{document}

%% file: macros.tex
\newcommand{\gl}{\lambda}
\newcommand{\go}{\omega}

\newcommand{\gO}{\Omega}

\newcommand{\cA}{\mathcal{A}}

\newcommand{\cE}{\mathcal{E}}

\newcommand{\cN}{\mathcal{N}}

\newcommand{\R}{\mathbb{R}}

\newcommand{\OAP}{$(\Omega,\mathcal{A},P)$}

\newcommand{\Ind}{\mathds{1}}

\newcommand{\lano}{\mathrm{o}}

\newcommand{\WIP}{\operatorname{WIP}}